\documentclass[12pt]{amsart}
\usepackage{amsmath,amscd,amssymb,amsfonts,graphics,mathrsfs}
\usepackage{hyperref}
\newcommand{\hl}{\hyperlink}
\newcommand{\htt}{\hypertarget}
\newcommand{\h}{\hbox}
\newcommand{\q}{\quad}
\newcommand{\nin}{\noindent}
\newcommand{\bs}{\par\bigskip}
\newcommand{\ms}{\par\medskip}
\newcommand{\sk}{\par\smallskip}
\newcommand{\bsn}{\par\bigskip\noindent}
\newcommand{\msn}{\par\medskip\noindent}
\newcommand{\skn}{\par\smallskip\noindent}
\newcommand{\ges}{\geqslant}
\newcommand{\les}{\leqslant}
\newcommand{\1}{\hskip1pt}

\newcommand{\msum}{\hbox{$\sum$}}
\newcommand{\mprod}{\hbox{$\prod$}}
\newcommand{\B}{{\mathscr B}}
\newcommand{\D}{{\mathscr D}}

\newcommand{\Hc}{{\mathscr H}}

\newcommand{\Lc}{{\mathscr L}}
\newcommand{\OO}{{\mathscr O}}
\newcommand{\Bf}{{\mathscr B}_{\!f}}
\newcommand{\Bt}{\widetilde{\mathscr B}}
\newcommand{\Btf}{\widetilde{\mathscr B}_{\!f}}
\newcommand{\Cf}{{\mathscr C}^{\ssb}_{\!f}}
\newcommand{\Ctf}{\widetilde{\mathscr C}^{\ssb}_{\!f}}
\newcommand{\Ctfo}{\widetilde{\mathscr C}^{\ssb}_{\!f_1}}
\newcommand{\Gf}{{\mathscr G}_{\!\1f}}

\newcommand{\Hsf}{{\mathscr H}''_{\!f}}

\newcommand{\Vt}{\widetilde{V}}
\newcommand{\Yt}{\widetilde{Y}}
\newcommand{\alt}{\widetilde{\alpha}}

\newcommand{\PP}{{\mathbb P}}
\newcommand{\Q}{{\mathbb Q}}
\newcommand{\C}{{\mathbb C}}
\newcommand{\N}{{\mathbb N}}

\newcommand{\Z}{{\mathbb Z}}

\newcommand{\Gr}{{\rm Gr}}
\newcommand{\DR}{{\rm DR}}

\newcommand{\al}{\alpha}
\newcommand{\be}{\beta}

\newcommand{\la}{\lambda}

\newcommand{\Om}{\Omega}

\newcommand{\dd}{\partial}
\newcommand{\ddd}{{\rm d}}
\newcommand{\eq}{\,{=}\,}
\newcommand{\gess}{\,{\ges}\,}
\newcommand{\less}{\,{\les}\,}
\newcommand{\sgt}{\,{>}\,}
\newcommand{\slt}{\,{<}\,}
\newcommand{\nes}{\,{\ne}\,}
\newcommand{\notins}{\,{\notin}\,}
\newcommand{\mi}{\1{-}\1}

\newcommand{\pl}{\1{+}\1}

\newcommand{\bl}{\bigl}
\newcommand{\br}{\bigr}
\newcommand{\sst}{\,{\subset}\,}

\newcommand{\ins}{\,{\in}\,}
\newcommand{\tos}{\,{\to}\,}
\newcommand{\defs}{\,{:=}\,}
\newcommand{\ssc}{\,\raise.15ex\hbox{${\scriptstyle\circ}$}\,}
\newcommand{\ssb}{\raise.15ex\h{${\scriptscriptstyle\bullet}$}}
\newcommand{\into}{\hookrightarrow}
\newcommand{\onto}{\twoheadrightarrow}
\newcommand{\simto}{\,\,\rlap{\hskip1.3mm\raise1.4mm\hbox{$\sim$}}\hbox{$\longrightarrow$}\,\,}
\begin{document}
\h{}\bs
\centerline{\large Examples of Hirzebruch-Milnor classes of projective hypersurfaces}
\sk
\centerline{\large detecting higher du Bois or rational singularities}
\bs
\centerline{Morihiko Saito}
\bs\bs\vbox{\nin\narrower\smaller
{\bf Abstract.} We show that it is possible to utilize the Hirzebruch-Milnor classes of projective hypersurfaces in the classical sense to detect higher du Bois or rational singularities only in some special cases. We also give several remarks clarifying some points in my earlier papers.}
\bs\ms
\centerline{\bf Introduction}
\bsn
The notion of {\it Hirzebruch-Milnor characteristic class\1} for a hypersurface $X$ in a smooth variety $Y$ was introduced in \cite{MSS1} applying the vanishing cycle functor of mixed Hodge modules (\cite{mhp}, \cite{mhm}) to a {\it global\1} defining function of $X$ in $Y$. Recently this class was employed in \cite{MY} to detect {\it higher du Bois\1} or {\it rational\1} singularities (\cite{JKSY3}, \cite{FL}) using in an essential way their characterizations in terms of {\it minimal exponents,} see \cite{MOPW}, \cite{JKSY3}, \cite{MP}, \cite[Appendix]{FL} (and also \cite{rat}, \cite{mos} for the case $k\eq0$). It is, however, quite non-trivial whether we have its {\it good examples,} since it is not easy to construct an example of a smooth complex variety $Y$ with a morphism $f\,{:}\,Y\tos\C$ such that $X\defs f^{-1}(0)$ is {\it analytic-locally irreducible\1} and its singular locus is {\it non-isolated\1} and {\it projective.} Recall that if $X$ has only {\it rational\1} singularities, then $X$ is {\it normal,} hence {\it analytic-locally irreducible.} Note also that higher du Bois (and rational) hypersurface singularities are rational singularities, and that the spectral Hirzebruch-Milnor class is equivalent to the {\it spectrum\1} in the {\it isolated\1} singularity case.
\sk
In these notes, we show that there are nevertheless some non-trivial examples (rather surprisingly) even though these are quite artificial and are restricted to entirely special situations, see Theorem~\hl{T1.3}{1.3} below. Note that the Hirzebruch-Milnor classes in the {\it new\1} sense are {\it directly\1} applicable to the detection of higher du Bois or higher rational singularities for {\it any\1} hypersurfaces of smooth projective varieties without taking any blow-ups nor assuming a condition on $k$ as in Thm.\,\hl{T1.3}{1.3} below, see \cite[Thm.\,1--2]{MSY}.
\bs\bs
\centerline{\bf 1. Smoothing of hypersurfaces}
\par\htt{1.1}{}\bsn
{\bf 1.1.~Deformation of hypersurfaces.} Firstly we have the following.
\par\htt{L1.1}{}\msn
{\bf Lemma~1.1.} {\it Let $Y$ be a smooth projective variety with $\Lc$ a very ample invertible sheaf. For $f,g\in\Gamma(Y,\Lc)$, consider a one-parameter deformation
$$V_{f+sg}:=\{f\pl sg\eq0\}\sst Y{\times}\C,$$
\par\nin where $s$ is the coordinate of $\C$. Then the total space $V_{f+sg}$ is singular if $V_f\defs\{f\eq0\}\sst Y$ has {\it non-isolated\1} singularities.}
\ms
This follows from the equality
\htt{1.1.1}{}
$${\rm Sing}\,V_{f+sg}\cap\{s\eq0\}=V_g\cap{\rm Sing}\, V_f.
\leqno(1.1.1)$$
\par\nin \par\htt{1.2}{}\msn
{\bf 1.2.~Desingularization of total spaces.}  We then observe the following.
\par\htt{L1.2}{}\msn
{\bf Lemma~1.2.} {\it If we desingularize the total space $V_{f+sg}$ in Lemma~{\rm\hl{L1.1}{1.1}} using blow-ups, then the total transform of $V_f\eq V_{f+sg}\cap\{s\eq0\}\sst V_{f+sg}$ has non-rational singularities.}
\ms
Indeed, the {\it exceptional divisors\1} of blow-ups must be contained in the {\it total transform of the central fiber,} hence the latter must be {\it reducible.} Here we assume the following two conditions:
\par\htt{C1}{}\msn
\vbox{\nin\rlap{(C1)}\hskip1cm\hangindent=1cm\hangafter=1
General members of the one-parameter family are smooth.
\par\htt{C2}{}\skn\rlap{(C2)}\hskip1cm\hangindent=1cm\hangafter=1
The desingularization induces an isomorphism over the smooth part.}
\sk
These imply that
\par\htt{C3}{}\skn\rlap{(C3)}\hskip1cm\hangindent=1cm\hangafter=1
The {\it centers\1} and {\it exceptional divisors\1} of blow-ups are contained in the central fibers.
\sk
Note that condition~(\hl{C1}{C1}) is satisfied if $g\ins\Gamma(Y,\Lc)$ is sufficiently general.
\par\htt{R1.2}{}\msn
{\bf Remark~1.2.} If we consider a morphism $h:\C^n\to\C$ defined by a polynomial denoted also by $h$, and take a natural compactification, this is similar to the above situation, since it is the one-parameter deformation defined by $f\mi s x_{n+1}^d$ ($s\ins\C)$ with $f$ a homogeneous polynomial of $x_1,\dots,x_{n+1}$ such that $f|_{x_{n+1}=1}\eq h$ and $d\defs\deg h=\deg f$.
\par\htt{1.3}{}\msn
{\bf 1.3.~Desingularization with condition~{\rm(\hl{C3}{C3})} unsatisfied.} The situation is, however, quite different from \hl{1.2}{1.2} {\it if condition~{\rm(\hl{C3}{C3})} is unsatisfied for all the centers and exceptional divisors of blow-ups.} It seems rather difficult to find a good example with condition~(\hl{C2}{C2}) unsatisfied (or $\dim{\rm Sing}\,X\gess 2$), since the centers of blow-ups should not be discrete and we have to construct them {\it globally.} As for condition (\hl{C1}{C1}), it is not necessarily easy to construct an example where the singularities of the total space can be resolved by blow-ups whose centers are {\it dominant\1} over the base space of the one-parameter family. We can nevertheless prove the following.
\par\htt{P1.3}{}\msn\vbox{\nin
{\bf Proposition~1.3.} {\it In the situation of Lemma~{\rm\hl{L1.1}{1.1},} assume the singular locus ${\rm Sing}\,V_f$ is a curve, and transversal slices of $V_f$ at general singular points are ordinary $m$-ple points with $m\gess 2$. Assume further that there is an ample invertible sheaf $\Lc'$ on $Y$ together with a section $h\ins\Gamma(Y,\Lc')$ such that $\Lc'{}^{\otimes m}\cong\Lc$, $V_h\sst Y$ is reduced and smooth, and intersects transversally ${\rm Sing}\,V_f$ at general points. Set $g\defs h^m$. Consider the blow-up $\Yt$ of $Y$ with reduced center $V_h\cap{\rm Sing}\,V_f\sst Y$. Let $\Vt_f$, $\Vt_{f+sg}$ be the proper transforms of $V_f$ and $V_{f+sg}$ in $\Yt$ and $\Yt{\times}\C$ respectively. Then $\Vt_{f+sg}$ is smooth, $\Vt_f\eq\Vt_{f+sg}\cap\{s\eq0\}$, and transversal slices of $\Vt_f$ at general singular points are ordinary $m$-ple points.}}
\msn
{\it Proof.} There are local coordinates $x_1,\dots,x_n$ of $Y$ ($n\defs\dim Y$) such that we have locally
\htt{1.3.1}{}
$$f\pl sg\eq\msum_{j\ges0}\,f_{j+m}(x_1,\dots,x_{n-1},x_n)\pl s\1x_n^m,
\leqno(1.3.1)$$
\par\nin trivializing the invertible sheaf $\Lc$, where the $f_j$ are homogeneous polynomials in $x_1,\dots,x_{n-1}$ of degree $j$ with coefficients in $\C\{x_n\}$. After the point center blow-up $\pi$ with $\pi^*x_i\eq y_iy_n$ ($i\slt n$), $\pi^*x_n\eq y_n$, we get locally
\htt{1.3.2}{}
$$\pi^*(f\pl sg)\eq\bl(\msum_{j\ges0}\,f_{j+m}(y_1,\dots,y_{n-1},y_n)y_n^j\pl s\br)y_n^m.
\leqno(1.3.2)$$
\par\nin So the assertion holds.
\ms
In Proposition~\hl{P1.3}{1.3}, conditions (\hl{C1}{C1}) and (\hl{C3}{C3}) are unsatisfied. The minimal exponents of $V_f$ and $\Vt_f$ {\it at general singular points\1} are equal to $\tfrac{n-1}{m}$ as a consequence of \cite{DMST} and \cite{exp}. Since the singularities do not change outside the center and exceptional divisor of blow-up, we then get the following.
\par\htt{T1.3}{}\msn\vbox{\nin
{\bf Theorem~1.3.} {\it With the notation and assumptions of Proposition~{\rm\hl{P1.3}{1.3}}, we can decide whether $V_f$ has only $k$-du Bois singularities by employing the Hirzebruch-Milnor class associated with the divisor $\Vt_f\eq\Vt_{f+sg}\cap\{s\eq0\}$ as in \cite{MY} if $k\les\tfrac{n-1}{m}\mi1$. A similar assertion holds for $k$-rational singularities if $k\slt\tfrac{n-1}{m}\mi1$.}}
\par\htt{R1.3a}{}\msn
{\bf Remark~1.3a.} Theorem~\hl{T1.3}{1.3} is not good for {\it explicit\1} calculations, since the definition of Hirzebruch-Milnor classes \cite{MSS1} is rather involved. It is better to calculate the {\it minimal exponent\1} $\alt_X$, which can be defined up to sign as the maximal root of the {\it reduced Bernstein-Sato polynomial\1} $b_f(s)/(s{+}1)$ for local defining functions $f$ of $X\sst Y$. Recall that
\par\htt{1.3.3}{}\skn\rlap{(1.3.3)}\hskip1.5cm\hangindent=1.5cm\hangafter=1
$X$ has only $k$-du~Bois (resp. $k$-rational) singularities if and only if $\alt_X\gess k{+}1$ (resp.\ $\alt_X\sgt k{+}1$),
\msn
see \cite{MOPW}, \cite{JKSY3}, \cite{MP}, \cite[Appendix]{FL} (and also \cite{rat}, \cite{mos} for the case $k\eq0$). This is closely related to the {\it microlocal $V$\!-filtration,} see \cite{hi}.
\par\htt{R1.3b}{}\msn
{\bf Remark~1.3b.} It is not necessarily easy to construct a concrete example for Theorem~\hl{T1.3}{1.3}. The latter theorem can be generalized to the case $\Lc'^{\otimes m'}\cong\Lc$ and $g\eq h^{m'}$ with $r\defs\tfrac{m'}{m}\ins\Z_{\ges2}$, where we have to repeat the point center blow-ups $r$ times in order to resolve the singularities. As an example for this generalized version, we can consider a divisor $D\sst\PP^{n-1}$ having only isolated singularities which are ordinary $m$-ple points. Let $f\in\C[x_1,\dots,x_n]$ be a defining polynomial of $D$. Assume $m'\eq d\defs\deg f\ins m\1\Z_{\ges2}$. Set $Y\defs\PP^n$ with $h\defs x_{n+1}$. Then we can apply Theorem~\hl{T1.3}{1.3} to this. In most cases the minimal exponent of $f$ is equal to $\tfrac{n}{d}\,\bl({<}\,\tfrac{n-1}{m}\br)$ assuming $n\gess 3$. As a concrete example, we may consider for instance 
\htt{1.3.4}{}
$$f\eq\msum_{k\les a}\,\msum_{i\ne k}\,x_i^mx_k^{d-m}+\msum_{k>a}\,x_k^d\q(a\ins[1,n]).
\leqno(1.3.4)$$
\par\nin \sk
For the example (\hl{1.3.4}{1.3.4}), it is expected that the minimal exponent $\alt_X$ is equal to $\tfrac{n}{d}$ using the theory of {\it pole order filtrations\1} as in \cite{DS1}, \cite{wh}, \cite{DS2}, \cite{nwh}. Here {\it the pole order spectral sequence degenerates at\1} $E_2$, and the computation becomes quite simple (see \cite{wh}), since the associated projective hypersurface has only {\it homogeneous\1} isolated singularities. Indeed, for $j\less a$, we have locally
\htt{1.3.5}{}
$$f|_{x_j=1}\eq\msum_{i\ne j}\,u_ix_i^m\eq\msum_{i\ne j}\,(u_i^{1/m}x_i)^m,
\leqno(1.3.5)$$
\par\nin with $u_i$ unit. In the case $m\eq2$, we can apply \cite[Thm.\,3]{DS2}.
\par\htt{R1.3c}{}\msn
{\bf Remark~1.3c.} It seems quite difficult to generalize Proposition~\hl{P1.3}{1.3} to the case where the transversal slices at general singular points are {\it not ordinary $m$-ple points,} even to the {\it weighted homogeneous\1} case. Here we cannot employ {\it weighted blow-ups,} since Hodge modules do not work well on $V$\!-manifolds. If we apply a (reduced) point center blow-up, then the minimal exponent may {\it decrease\1} in general.
\sk
Consider for instance the case we have locally
\htt{1.3.6}{}
$$f\eq\msum_{i=1}^{n-2}\,x_i^a+x_{n-1}^b,\q g\eq x_n^b\q\h{with}\q a\sgt b\sgt 1,
\leqno(1.3.6)$$
\par\nin trivializing the invertible sheaf $\Lc$. After the point center blow-up $\pi$ with $\pi^*x_i\eq y_iy_n$ ($i\slt n$), $\pi^*x_n\eq y_n$, we get locally
\htt{1.3.7}{}
$$\pi^*(f{+}sg)\eq\bl(\bl(\msum_{i=1}^{n-2}\,y_i^a\br)\1y_n^{a-b}\pl y_{n-1}^b\pl s\br)y_n^b.
\leqno(1.3.7)$$
\par\nin Note that the minimal exponents of $h_1\defs\msum_{i=1}^{n-2}\,y_i^a$ and $h_2\defs y_n^{a-b}$ are respectively given by
\htt{1.3.8}{}
$$\alt_{h_1}\eq\tfrac{n-2}{a},\q\alt_{h_2}\eq\tfrac{1}{a-b}.
\leqno(1.3.8)$$
\par\nin We can see that $\alt_h\eq\min\{\alt_{h_1},\alt_{h_2}\}$ with $h\defs h_1h_2$. So the minimal exponent {\it decreases\1} if
\htt{1.3.9}{}
$$\tfrac{n-2}{a}\sgt\tfrac{1}{a-b}\q\h{or equivalently}\q\tfrac{a}{b}\sgt\tfrac{n-2}{n-3}.
\leqno(1.3.9)$$
\par\nin Here we apply the Thom-Sebastiani type theorem for Bernstein-Sato polynomials as in \cite[Thm.\,0.8]{mic} (where one of the functions must be weighted homogeneous).
\bs\bs
\centerline{\bf 2. Some remarks}
\par\htt{2.1}{}\bsn
{\bf 2.1.~Microlocal $V\!$-filtration.} The assertion in \cite[Prop.\,3.14]{MY} immediately follows from the canonical isomorphisms shown in \cite[(2.2.1) and Lem.\,2.2]{mic}\,:
\htt{2.1.1}{}\par\vbox{
$$\!\!\!\!\!\!\!\!\!\!\Gr_V^{\al}\iota:\Gr_V^{\al}(\Bf,F)\simto\Gr_V^{\al}(\Btf,F)\q\q(\forall\,\al<1),
\leqno(2.1.1)$$
\par\nin \vskip-6mm
\htt{2.1.2}{}
$$\q\q\q\q\q\!\dd_t^j:F_pV^{\al}\Btf\simto F_{p+j}V^{\al-j}\Btf\q\q\!(\forall\,j,p\ins\Z,\,\al\ins\Q).
\leqno(2.1.2)$$}
Here $\iota:\Bf\defs\OO_Y[\dd_t]\delta(t{-}f)\into\Btf\defs\OO_Y[\dd_t,\dd_t^{-1}]\delta(t{-}f)$ is the canonical inclusion, $F$ is the filtration by the order of $\dd_t$, and $V$ is the (microlocal) $V$\!-filtration of Kashiwara and Malgrange.
\sk
As for the proof of (\hl{2.1.1}{2.1.1}), note that
\htt{2.1.3}{}
$$F_p\Gr_V^{\al}\Btf\eq F_pV^{\al}\Btf/F_pV^{>\al}\Btf\q\q(\forall\,p\ins\Z,\,\al\ins\Q),
\leqno(2.1.3)$$
\par\nin (similarly for $\Bf$), and we have the canonical isomorphisms
\htt{2.1.4}{}
$$\iota:F_pV^{\al}\Bf\simto F_pV^{\al}\Btf/F_{-1}\Btf\q\q(\forall\,p\gess 0,\,\al\less1).
\leqno(2.1.4)$$
\par\nin Indeed, $F_{-1}\sst F_p$, and we can easily deduce the {\it equalities\1} as submodules of $\Btf$
\htt{2.1.5}{}
$$F_pV^{\al}\Btf\eq F_{-1}\Btf\oplus\iota\bl(F_pV^{\al}\Bf\br)\q\q(\forall\,p\gess 0,\,\al\less1)
\leqno(2.1.5)$$
\par\nin from the definition of $V$ on $\Btf$ (see \cite[(2.1.3)]{mic}), that is,
\htt{2.1.6}{}
$$V^{\al}\Btf\eq F_{-1}\Btf\oplus\iota\bl(V^{\al}\Bf\br)\q\q(\forall\,\al\less 1),
\leqno(2.1.6)$$
\par\nin taking the intersection with $F_p\Btf$ and using the {\it strict compatibility\1} of $\iota$ with $F$ (and also $(A{+}B)\,{\cap}\,C\eq A\pl B\,{\cap}\,C$ if $A\sst C$).
\sk
As for \cite[Lem.\,3.13]{MY}, which asserts in view of (\hl{2.1.2}{2.1.2}) the canonical isomorphisms
\htt{2.1.7}{}
$$\Gr^F_p\iota:\Gr_p^FV^{\al}\Bf\simto\Gr_p^FV^{\al}\Btf\q\q(\forall\,p\gess0,\,\al\slt1),
\leqno(2.1.7)$$
\par\nin this also follows easily from (\hl{2.1.5}{2.1.5}), since the latter or (\hl{2.1.6}{2.1.6}) implies that $F_{-1}V^{\al}\Btf\eq F_{-1}\Btf$ for $\al\less 1$, taking the intersection with $F_{-1}\Btf$.
\par\htt{2.2}{}\msn
{\bf 2.2.~Explicit description of the microlocal $V$\!-filtration.} In the isolated singularity case the $V$\!-filtration on the Gauss-Manin system $\Gf$ can be induced from the $V$\!-filtration on $\Bf$ taking the relative de Rham complex (see for instance \cite[3.4.8]{mhp}), where we can replace $\Bf$ with $\Btf$ as in \cite[(4.11.5)]{DS1}. Indeed, we have the isomorphisms
$$\aligned\Cf/V^{\alpha}\Cf&\simto\Ctf/V^{\alpha}\Ctf\q\h{with}\q H^{-1}\Cf/V^{\alpha}\Cf\eq0\q(\forall\,\al\less1),\\ \dd_t^j\,{:}\,V^{\al}\Gf&\simto V^{\al-j}\Gf\q(\forall\,\al\ins\Q,\,j\ins\Z).\endaligned$$
\par\nin Here $\Cf\defs\DR_{Y\times\C/\C}(\B_{\!f,0})$ denotes the {\it Gauss-Manin complex,} $\Ctf\defs\DR_{Y\times\C/\C}(\Bt_{\!f,0})$ is the algebraic microlocal one, and the {\it Gauss-Manin system\1} is defined by $\Gf\defs H^0\Cf\simto H^0\Ctf$. However, the $V$\!-filtration on the {\it Brieskorn lattice}
$$\Hsf\defs\Om_{Y,0}^n/\ddd f{\wedge}\1\ddd\Om^{n-2}_{Y,0}\q(n\defs\dim Y)$$
\par\nin is {\it never\1} induced from the {\it microlocal\1} $V$\!-filtration $\Vt$ on $\OO_Y$ trivializing $\Om_Y^n$ (unless $\Gr^F_0$ is taken with $F$ the {\it microlocal\1} Hodge filtration), and the filtration $\Vt$ cannot be captured by using $b_{f,g}(s)/(s{+}1)$ as in \cite{BBV} assuming $g\notins f\OO_Y$.
\sk
Consider for instance the case $f\eq\msum_{i=1}^n\,x_i^{m_i}$, $g\eq\mprod_{i=1}^n\,x_i^{a_i}$ ($m_i\gess 2$) with $a_i{+}1\notins m_i\1\Z$ for any $i\ins[1,n]$, that is, $[g\ddd x]\nes0$ in $\Hsf$. Here $n\eq d_Y$, $\ddd x\defs\ddd x_1{\wedge}\cdots{\wedge}\ddd x_n$. Let $\alt_f(g)$ be the maximal root of $b_{f,g}(s)/(s{+}1)$ up to sign. Let $\al^{\Vt}_f(g)$, $\al^{\rm Br}_f(g)$ be respectively the maximum of $\al\ins\Q$ with $g\ins\Vt^{\al}\OO_Y$ and $[g\ddd x]\ins V^{\al}\Hsf$. Let $b_i,c_i\in\N$ such that $a_i\eq b_i\pl c_i\1(m_i{-}1)$ with $b_i\ins[0,m_i{-}2]$. We can verify that
\htt{2.2.1}{}
$$\al^{\Vt}_f(g)\eq\al^{\Vt}_f(x^b)\pl\msum_ic_i\eq\msum_i\bl(\tfrac{b_i+1}{m_i}\pl c_i\br),\q\al^{\rm Br}_f(g)\eq\msum_i\tfrac{a_i+1}{m_i},
\leqno(2.2.1)$$
\par\nin where the first equality follows from the Thom-Sebastiani type theorem for microlocal $V$\!-filtrations \cite{MSS2}, see also \cite[(1.1.7)]{JKSY2}, \cite{hi}. We thus get that
\htt{2.2.2}{}
$$\al^{\Vt}_f(g)=\al^{\rm Br}_f(g)\pl\msum_i\,\tfrac{c_i}{m_i}.
\leqno(2.2.2)$$
\par\nin \sk
Assume there is a relation
\htt{2.2.3}{}
$$b(s)gf^s\eq Pfgf^s\q\h{with}\q b(s)\ins\Q[s],\,P\ins\D_Y[s].
\leqno(2.2.3)$$
\par\nin This implies a relation in the Gauss-Manin system $\Gf\eq \Hsf[\dd_t]$ (which is the localization of $\Hsf$ by $\dd_t^{-1}$). Considering this relation modulo $V^{>\al}\Gf$ with $\al\defs\al^{\rm Br}_f(g)\eq\msum_i\tfrac{a_i+1}{m_i}$, we see that
\htt{2.2.4}{}
$$\h{$-\al^{\rm Br}_f(g)$ is a root of $b(s)$,}
\leqno(2.2.4)$$
\par\nin since the right-hand side of the induced relation is equal to $[P_0fg\ddd x]$, and vanishes, where $P_0\ins\OO_Y[s]$ coincides with $P$ modulo $\msum_i\,\dd_{x_i}\D_Y[s]$. Recall that $\Gf$ is the top cohomology of the relative de Rham complex of $\Bf\eq\OO_Y[\dd_t]\delta(t{-}f)$ for $Y{\times}\C\to\C$, where $f^s$ is identified with $\delta(t{-}f)$, and $s$ with $-\dd_tt$.
\sk
By the definition of $\alt_f(g)$ we then get the inequality
\htt{2.2.5}{}
$$\alt_f(g)\less\al^{\rm Br}_f(g).
\leqno(2.2.5)$$
\par\nin Combined with (\hl{2.2.2}{2.2.2}), this implies the inequality
\htt{2.2.6}{}
$$\al^{\Vt}_f(g)\mi\alt_f(g)\gess\msum_i\,\tfrac{c_i}{m_i}.
\leqno(2.2.6)$$
\par\nin \par\htt{R2.2a}{}\msn
{\bf Remark~2.2a.} By \cite[Thm.\,2]{JKSY2} we have in the general singularity case
\htt{2.2.7}{}
$$\Vt^{\al}\OO_Y\sst(\dd f)\,\,\,\,\h{if}\,\,\,\,\al\sgt d_Y\mi\alt_f,
\leqno(2.2.7)$$
\par\nin where $(\dd f)$ denotes the Jacobian ideal of $f$. This implies that
\htt{2.2.8}{}
$$f^kg\ins(\dd f)\,\,\,\,\h{if}\,\,\,\,k\sgt d_Y\mi\alt_f\mi\al^{\Vt}_f(g).
\leqno(2.2.8)$$
\par\nin Note that $\al^{\Vt}_f(g)$ can be replaced by the maximum of $\al^{\Vt}_f(h)$ with $h\mi g\ins(\dd f)$. In other words it is enough to determine the {\it quotient\1} $V$\!-filtration on $\C\{x\}/(\dd f)$. In the {\it isolated\1} singularity case, this quotient filtration coincides with the one coming from the $V$\!-filtration on the Brieskorn lattice $\Hsf$, see \cite[4.11]{DS1}, \cite{JKSY1}.
\par\htt{R2.2b}{}\msn
{\bf Remark~2.2b.} The definition of microlocal $V$\!-filtration is misstated just after \cite[Prop.\,3.2.1]{BBV}, where the {\it intersection\1} is used as in the first version of \cite{MY} instead of the {\it graded quotient\1} $\Gr^F_0$ of the {\it microlocal\1} Hodge filtration $F$ on $\Btf$ (although \cite{BBV} was not quoted there).
\par\htt{R2.2c}{}\msn
{\bf Remark~2.2c.} The induced filtration $V$ by the inclusion $\OO_Y\eq F_0\Bf\into\Bf$ coincides with the {\it multiplier ideals\1} having the periodicity by multiplication by $f$ (see \cite{BS}), and does {\it not\1} induce the $V$\!-filtration on $\Hsf$ nor on $\C\{x\}/(\dd f)$ in the isolated singularity case. We have to use the induced filtration $V$ by the inclusions $F_p\Bf\sst\Bf$ for $p\gess 0$ in order to get the filtration $V$ on $\Hsf$ and on $\C\{x\}/(\dd f)$. Note that the bifiltered Gauss-Manin complex
$$(\Cf;F,V)\eq\DR_{Y\times\C/\C}(\Bf;F,V)$$
\par\nin is {\it never bistrict\1} even in the isolated singularity case, and the {\it spectral numbers\1} contained in $(p\1{-}1,p\1]$ are determined by the $\dim\Gr_F^{n-p}H^{n-1}(F_{\!f},\C)_{\la}$ ($\la\ins\C^*$) with $F_{\!f}$ the Milnor fiber, see \cite{St}, \cite{ovi}.
\par\htt{R2.2d}{}\msn
{\bf Remark~2.2d.} In the semi-weighted-homogeneous isolated singularity case, it is known that the $V$\!-filtration on the Gauss-Manin system $\Gf$ is induced from the (shifted) Newton filtration $V_w$ defined by using the weights $w_i$, see \cite{exp}, \cite{JKSY1}. However, this does not seem to imply immediately that the induced $V$\!-filtration on the Brieskorn lattice $\Hsf$ is the {\it quotient\1} filtration of the shifted Newton filtration $V_w$ on $\OO_Y$ (trivializing $\Om_Y^n$), and some more argument seems to be needed. Indeed, the assertion for the Gauss-Manin system implies the {\it inclusion\1} $V_w^{\ssb}\sst V^{\ssb}$ on $\Hsf$, so it is enough to show the equality of the dimensions of their graded pieces which are closely related to the multiplicities of spectral numbers. Using the algebraic microlocal Gauss-Manin complex $\Ctf$ and the invariance of the spectrum under $\mu$-constant deformations, this equality can be reduced to the case of the associated weighted homogeneous polynomial $f_1$ (where the assertion is easily shown using the stability of $\Hc''_{\!f_1}$ by the action of $\dd_tt$ via the interior product with the Euler field as in \cite[1.4]{BaSa}, see also \cite[1.5]{JKSY1}). Note that the graded pieces of $V_w$ on the algebraic microlocal Gauss-Manin complex $\Ctf$ are canonically isomorphic to those of $\Ctfo$ in a compatible way with the microlocal Hodge filtration $F$, and these microlocal complexes endowed with the two filtrations $F,V_w$ are {\it bistrict,} that is, {\it strict\1} in the sense of \cite[1.2.1]{mhp} (hence we have the commutativity of $H^0$, $\Gr^F_p$, $\Gr_V^{\al}$, see \cite[1.2.13]{mhp}) by an argument similar to \cite{exp}, see also \cite[4.11]{DS1}. Here we use the commutative diagram of exact sequences for any $p\slt q$, $\al\sgt\be$\,{:}
$$\begin{array}{cccccccccccccccc}H^0F_pV_w^{\al}\Ctf&\into&H^0F_pV_w^{\be}\Ctf&\onto&H^0F_p(V_w^{\be}/V_w^{\al})\Ctf\\ \rotatebox{-90}{\!\!\!\!$\into$\,}&&\rotatebox{-90}{\!\!\!\!$\into$\,}&&\rotatebox{-90}{\!\!\!\!$\into$\,}\\ H^0F_qV_w^{\al}\Ctf&\into&H^0F_qV_w^{\be}\Ctf&\onto&H^0F_q(V_w^{\be}/V_w^{\al})\Ctf\\
\end{array}$$
\par\nin where the injectivity of the vertical and left horizontal morphisms follows from
$$H^{-1}(F_q/F_p)(V_w^{\be}/V_w^{\al})\Ctf\eq H^{-1}(F_q/F_p)V_w^{\al}\Ctf\eq H^{-1}F_p(V_w^{\be}/V_w^{\al})\Ctf\eq0.$$
\par\nin Note that the graded pieces of the microlocal Hodge filtration $F$ on $\Ctf$ are the Koszul complex for the partial derivatives $\dd_{x_i}f$ which form a regular sequence.
\par\htt{2.3}{}\msn
{\bf 2.3.~Positivity of Todd class transformation.} The meaning of ``the positivity of the Todd class transformation on a projective space" in the proof of Theorem~A in \cite[\S4.2]{MY} does not seem very clear to the reader. It does not seem to be recommended to formulate this ``positivity" without using the notion of {\it topological filtration\1} on the Grothendieck group of coherent sheaves as in \cite{SGA6}, \cite{Fu}. One should note at least that if the support of a coherent sheaf on projective space has dimension $k$, then the degree $j$ part of its image by the Todd class transformation vanishes for $j\sgt 2k$ (by functoriality), and its degree $2k$ part is {\it positive,} see also \cite[\S1.3]{MSS3}, \cite[\S2.2]{MSY}. Note that this is a quite essential part of the proof of the main theorem.
\msn

\sk
{\smaller\smaller
RIMS Kyoto University, Kyoto 606-8502 Japan}
\end{document}